**WARING'S PROBLEM FOR FIELDS (*)**

**W.J. Ellison**

*To "Ian" Cassels on his 90th birthday*

**1. Introduction**. The classic Waring problem for $\boldsymbol{Z}$ has a vast literature (for a glimpse see [6], [18]). The solution of the corresponding problems associated with expressing the elements of a field $\boldsymbol{K}$ as the sum of $k^{\text{th}}$ powers is much less complete. For $k = 2$ we have Artin's characterisation that $\alpha \in \boldsymbol{K}$ is a sum of squares if and only if $\alpha > 0$ for all orderings of $\boldsymbol{K}$. For $k$ an even number, Becker [1] extended Artin's ideas to characterise $\alpha \in \boldsymbol{K}$ as a sum of $k^{th}$ powers if and only if $\alpha > 0$ for all orderings of $\boldsymbol{K}$ and $k \mid \upsilon(\alpha)$ for all valuations $\upsilon$ of $\boldsymbol{K}$ with a formally real valuation field.

Denote by $P(\boldsymbol{K}, k)$ the $a \in \boldsymbol{K}$ which are sums of $k^{\text{th}}$ powers of elements of $\boldsymbol{K}$, by $P^+(\boldsymbol{K}, k)$ the set of $a \in \boldsymbol{K}$ which are sums of $k^{\text{th}}$ powers of totally positive elements of $\boldsymbol{K}$. We are interested in deciding whether or not there exist integers $w(\boldsymbol{K}, k)$ and $g(\boldsymbol{K}, k)$ such that:

$a \in P(\boldsymbol{K}, k)$ implies that $a$ is the sum of at most $w(\boldsymbol{K}, k)$ $k^{\text{th}}$ powers;

$a \in P^+(\boldsymbol{K}, k)$ implies that $a$ is the sum of at most $g(\boldsymbol{K}, k)$ totally positive $k^{\text{th}}$ powers.

Neither Artin's nor Becker's characterisations gives any information about the existence of $w(k, \boldsymbol{K})$ or $g(K, k)$. This is to be expected since there are many fields for which $w(\boldsymbol{K}, 2) = \infty$, which implies that $w(\boldsymbol{K}, 2k) = \infty$.

The integer $w(\boldsymbol{K}, 2)$ is called the Pythagorean number of $\boldsymbol{K}$. For a given integer $s$, one can construct fields such that $w(\boldsymbol{K}, 2) = s$ [9]. However, for a given field it is usually a difficult problem to determine its Pythagorean number. For an algebraic number field $\boldsymbol{K}$ it is a classic result that $w(\boldsymbol{K}, 2) \leq 4$ and $w(\boldsymbol{K}(X), 2) \leq 5$ [16, 10]. For $n \geq 2$, $5+n \leq w(\boldsymbol{K}(X_1, \ldots, X_n), 2) \leq 2^{n+1}$ [2, 5]. If $\boldsymbol{R}$ is a real closed field, $w(\boldsymbol{R}(X), 2) = 2$, $w(\boldsymbol{R}(X_1, X_2), 2) = 4$ [3] and for $n \geq 3$, $1+ n \leq w(\boldsymbol{R}(X_1, \ldots, X_n), 2) \leq 2^n$ [2, 14, 15].

We will show that if $w(\boldsymbol{K}, 2)$ is finite and provided that the $k^{th}$ powers are dense (in a sense described explicitly in theorem 3) in $\boldsymbol{K}$, then $w(\boldsymbol{K}, k)$ is also finite for $k > 2$. The proofs are constructive, but the implied upper bounds for $w(\boldsymbol{K}, k)$ are large. This is to be expected since the method of proof does not use any deep arithmetical or algebraic properties of the field $\boldsymbol{K}$.

---

* Based upon part of a 1970 University of Cambridge Ph.D. dissertation [7], written under the direction of Pr. J.W.S. Cassels

**2. The basic theorems** Our first theorem treats the case when $k$ is odd. It is well known and is included here for completeness. When $k$ is even it simplifies the exposition if we separate the two cases: Theorem 2 (***K*** a non-real field i.e. $-1$ is a sum of squares in ***K***) and Theorem 3 (***K*** a formally real field), though the proofs are very similar. As an application of theorem 3 we give (Theorem 4) a characterisation of those rational functions over certain real fields ***K*** which can be written as sums of $k^{th}$ powers of rational functions. For these fields we show that $w(\mathbf{K}(X), 2) = g(\mathbf{K}(X), 2)$ (Theorem 5).

**Theorem 1.** *If **K** is a field of characteristic 0, then for every odd integer k,* $\mathbf{K} = P(\mathbf{K}, k)$ *and* $w(\mathbf{K}, k) < \infty$.

**Proof**. Take finite difference of $x^k$: If $\Delta(x^k) = (x+1)^k - x^k$; $\Delta^2(x^k) = \Delta(\Delta(x^k))$ etc., then by calculating successively the differences we have:

$$\Delta^k(x^k) = k!x + (k-1)k!/2 \text{ and } \Delta^k(x^k) = \sum_{r=0}^{k-1}(-1)^{k-1-r}\left(\frac{(k-1)....(k-r)}{r(r-1)...3.2.1}\right)(x+r)^k .$$

We thus have the identity over $\mathbf{Q}[x]$:

$$x = -\frac{(k-1)}{2} + \sum_{r=0}^{k-1}(-1)^{k-1-r}\frac{(k-1)(k-2)...(k-r)}{k!r(r-1)...1}(x+r)^k$$

Since $k$ is odd we can write the first term as a sum of at most $(k-1)/2$ $k^{\text{th}}$ powers and each of the rational numbers as the sum of $k^{\text{th}}$ powers.

**Theorem 2.** *If **K** is a non-real field of characteristic 0, then* $\mathbf{K} = P(\mathbf{K}, k)$ *and* $w(\mathbf{K}, k) < \infty$**.**

**Theorem 3**. *Let **K** be a formally real field. Suppose that* $\alpha \in \mathbf{K}$ *has the following two properties:*

*(i)* $\alpha$ *can be written as a sum of at most s squares in **K**;*

*(ii) there exists* $\beta \in \mathbf{K}$*, depending upon* $\alpha$*, that satisfies* $0 \prec \frac{s}{s+2k}\alpha \prec \beta^k \prec \alpha$ *for all orderings* $\rightleftarrows$ *of **K**,*

*then there exists a positive integer* $\gamma(s, k)$*, depending only on s and k such that:*

*(a)* $\alpha$ *can be written as a sum of* $\gamma(s, k)$ $k^{th}$ *powers of elements of **K**.*

*(b) If* $\beta$ *can be chosen to be totally positive, then* $\alpha$ *can be written as a sum of* $\gamma(s, k)$ $k^{th}$ *powers of totally positive elements of **K**.*

**Cor. 3.1.** *If* $\mathbf{K} = \mathbf{Q}$*, then* $w(\mathbf{Q}, k) \leq g(\mathbf{Q}, k) < \infty$ *for every exponent k.*
**Cor. 3.2.** *If **K** is an algebraic number field, then* $w(\mathbf{K}, k) \leq g(\mathbf{K}, k) < \infty$ *for every k*

The key result we need to prove theorems 2 and 3 is the existence of the 'Hilbert' identities:

**Lemma 1.** *For all positive integers k and s there is an integer* $M := M(k, s) = (2k+3)....(2k+2+s)/s!$*, positive rational numbers* $\lambda_I$ *with* $0 \leq i \leq M$ *and integers* $\alpha_{ij}$ *with* $0 \leq i \leq s$ *and* $0 \leq j \leq M$ *such that we have an identity of the form:*

(i) $\left(X_0^2+.....+X_s^2\right)^{k+1}=\sum_{j=1}^{M}\lambda_j\left(\alpha_{0j}X_0+....+\alpha_{sj}X_s\right)^{2k+2}$

(ii) $\left(X_0^2+....+X_s^2\right)^k+2kX_0^2(X_0^2+....+X_s^2)^{k-1}=(2k+1)\sum_{j=1}^{M}\lambda_j\alpha_{0j}^2\left(\alpha_{0j}X_0+....+\alpha_{sj}X_s\right)^{2k}$

**Proof of lemma 1.** For a short existence proof of (i), based upon the properties of convex sets, see [6]. For explicit constructive proofs see [12], [13] and [17]. The identity (ii) is obtained from (i) by differentiating twice with respect to $X_0$.

**Cor. 1**. *For any field* $\boldsymbol{K}$ *and positive integer k,* $P^+(\boldsymbol{K}, k) = P(\boldsymbol{K}, 2k)$.

**Proof of Cor. 1**. If $a \in P^+(\boldsymbol{K}, k)$, then $a = b_1^k + ... + b_n^k$, where the $b_i$ are totally positive. By Artin's theorem we have $b_i = c_1^2 + ... + c_r^2$ and we use the appropriate Hilbert identity to write $b_i^k$ as a sum of $2k$ th powers, so $a \in P(\boldsymbol{K}, 2k)$.

Conversely, if $a \in P(\boldsymbol{K}, 2k)$, then $a = b_1^{2k} + ... + b_s^{2k} = (b_1^2)^k + ... + (b_n^2)^k \in P^+(\boldsymbol{K}, k)$.

**Proof of theorem 2**. If $-1 = a_1^2 + ..... + a_s^2$ then any $\alpha \in \boldsymbol{K}$ can be written as :

$$\alpha=\left(\frac{\alpha+1}{2}\right)^2-\left(\frac{\alpha-1}{2}\right)^2=\left(\frac{\alpha+1}{2}\right)^2+\sum_{i=1}^{s}a_i^2\left(\frac{\alpha-1}{2}\right)^2$$

Thus, *every* element of $\boldsymbol{K}$ is a sum of at most $(s+1)$ squares. If $a_i \in \boldsymbol{K}$ then we have $1-a_i^2=b_{1i}^2+....+b_{(s+1)i}^2$, where $b_{ji} \in \boldsymbol{K}$ for $1 \leq j \leq (s+1)$. Substitute $X_0 = a_i$ and $X_j = b_{ji}$ for $1 \leq j \leq (s+1)$ in the identity (ii) to obtain:

$$1+2ka_i^2=(2k+1)\sum_{j=1}^{M}\lambda_j\alpha_{0j}^2(\alpha_{0j}a_i+\alpha_{1j}b_{1i}+....+\alpha_{sj}b_{si})^{2k}$$

If $a \in \boldsymbol{K}$ we have $a=a_1^2+.....+a_{(s+1)}^2$ with $a_i \in \boldsymbol{K}$. Substitute the $a_i$ into the above identity and sum the resulting set of $(s+1)$ equations to obtain:

$$s+1+2ka=\sum_{i=1}^{s+1}\sum_{j=1}^{M}(2k+1)\lambda_j\alpha_{0j}^2\left(\alpha_{0j}a_i+\alpha_{1j}b_{1i}+....+\alpha_{sj}b_{si}\right)^{2k}$$

Let $A_k$ be the least common multiple of the denominators of the $\lambda_j$ for $0 \leq j \leq M$, so that $A_k\lambda_j = \Lambda_j$, positive integer. We now have:

$$(s+1)A_k+2kA_ka=\sum_{i=1}^{s+1}\sum_{j=1}^{M}\Lambda_j\alpha_{0j}^2\left(2k+1\right)\left(\alpha_{0j}a_i+.....+\alpha_{sj}b_{si}\right)^{2k}=\sum_{i=1}^{n}\zeta_i^{2k},$$

where $\zeta_i \in \boldsymbol{K}$ for $1 \leq i \leq n$ and $n=(2k+1)(s+1)\sum_{j=1}^{M}\Lambda_j\alpha_{0j}^2$

For any $\alpha \in \boldsymbol{K}$ we take $a = (\alpha - (s+1)A_k)/2kA_k$ and substitute in the above equation we obtain: $\alpha=\sum_{i=1}^{n}\xi_i^{2k}$ and the assertion follows with $w(\boldsymbol{K}, k) = g(\boldsymbol{K}, k) \leq n$.

**Proof of theorem 3.** If $a_i \in \boldsymbol{K}$ and $1 - a_i^2 \nleq 0$ for all orderings '$\nleq$' of $\boldsymbol{K}$, then $1 - a_i^2$ is totally positive and by hypothesis (i) we have $1-a_i^2=b_{1i}^2+....+b_{si}^2$, where $b_{ji} \in \boldsymbol{K}$ for $1 \leq j \leq s$. Substitute $X_0 = a_i$ and $X_j = b_{ji}$ for $1 \leq j \leq s$ in the identity (ii) of lemma 1 to obtain:

$$1+2ka_i^2=(2k+1)\sum_{j=1}^{M}\lambda_j\alpha_{0j}^2(\alpha_{0j}a_i+\alpha_{1j}b_{1i}+....+\alpha_{sj}b_{si})^{2k}$$

If $a \in \boldsymbol{K}$ and $0 \rightleftarrows a \rightleftarrows 1$ for all orderings '$\rightleftarrows$' of $\boldsymbol{K}$, the by hypothesis (i) we have $a=a_1^2+.....+a_s^2$ with $a_i \in \boldsymbol{K}$. Substitute the $a_i$ into the above identity and sum the resulting set of $s$ equations to obtain:

$$s+2ka=\sum_{i=1}^{s}\sum_{j=1}^{M}(2k+1)\lambda_j\alpha_{0j}^2\left(\alpha_{0j}a_i+\alpha_{1j}b_{1i}+....+\alpha_{sj}b_{si}\right)^{2k}$$

Let $A_k$ be the least common multiple of the denominators of the $\lambda_j$ for $0 \leq j \leq M$, so that $A_k\lambda_j = \Lambda_j$, a positive integer. We now have:

$$sA_k+2kA_ka=\sum_{i=1}^{s}\sum_{j=1}^{M}\Lambda_j\alpha_{0j}^2\left(2k+1\right)\left(\alpha_{0j}a_i+.....+\alpha_{sj}b_{si}\right)^{2k}=\sum_{i=1}^{n}\zeta_i^{2k},$$

where $\zeta_i \in \boldsymbol{K}$ for $1 \le i \le n$ and $n = (2k+1)s\sum_{j=1}^{M} A_j \alpha_{0j}^2$

If $\alpha \in \boldsymbol{K}$ and $\alpha$ is totally positive, then so is $A_k(s+2k)/\alpha$. Hence, by hypothesis (ii) there is a $\beta \in \boldsymbol{K}$, which depends upon $\alpha$, such that

$$0 \prec \frac{s}{s+2k}\left(\frac{A_k(s+2k)}{\alpha}\right) \prec \beta^k \prec \left(\frac{A_k(s+2k)}{\alpha}\right)$$

for all orderings '$\prec$' of $\boldsymbol{K}$. If we put $a = (\alpha\beta^k - sA_k)/2kA_k$, then $0 \prec a \prec 1$ for all orderings '$\prec$' of $\boldsymbol{K}$. Substitute this $a$ in the above equation to obtain: $\alpha\beta^k = \sum_{i=1}^{n} \xi_i^{2k}$ and the assertions of the theorem follow with $\gamma(s, k) = n$.

**Proof of Cor. 3.1.** Any positive rational number is the sum of at most *4* squares of positive rational numbers. There is only one ordering on $\boldsymbol{Q}$ and for any positive integer $k$ and positive rational number $\alpha$ there is a positive rational number $\beta$ such that $2\alpha/(2+k) < \beta^k < \alpha$.

**Proof of Cor.3.2.** If $\boldsymbol{K}$ is a totally imaginary number field, then it is a classic theorem that every $\alpha \in \boldsymbol{K}$ is the sum of at most *4* squares and the result follows from Theorem 2. If $\boldsymbol{K}$ is an algebraic number field that is not totally imaginary, then it is a classic result that every totally positive element of $\boldsymbol{K}$ can be written as a sum of at most *4* squares. Thus, hypothesis (i) is satisfied.

We now show that (ii) can be satisfied with $\beta$ totally positive. If $\boldsymbol{K} = \boldsymbol{Q}(\theta)$ denote by $\boldsymbol{K}^{(1)} = Q(\theta^{(1)}), \ldots, \boldsymbol{K}^{(r)} = Q(\theta^{(r)})$ the real conjugate fields. We need the following trivial remark: "*If* $\varepsilon > 0$, $\eta_1, \ldots, \eta_r$ *be given real numbers, then there is a* $\beta \in \boldsymbol{K}$ *such that* $|\beta^{(i)} - \eta_i| < \varepsilon$ *for* $1 \le i \le r$.". The proof is simple: Let $f(x)$ be a polynomial of degree $(r-1)$ with real coefficients, which takes the values $\eta_i + \varepsilon/2$ at $x = \theta^{(i)}$ for $1 \le i \le r$. Let $g(x)$ be a polynomial of degree $(r-1)$ with rational coefficients such that that $|f(x) - g(x)| < \varepsilon/2$ for $x = \theta^{(1)}, \ldots, \theta^{(r)}$. Put $\beta = g(\theta)$, then $\beta^{(i)} = g(\theta^{(i)})$ for $1 \le i \le r$ and we have $|\beta^{(i)} - \eta_i| < \varepsilon$ for $1 \le i \le r$.

We now return to the verification of hypothesis (ii). If $\alpha \in \boldsymbol{K}$ and $\alpha$ is totally positive, then $\alpha^{(i)} > 0$ *for* $1 \le i \le r$. To satisfy hypothesis (ii) we must find a $\beta \in \boldsymbol{K}$, depending upon $\alpha$, such that $0 < \left(\frac{2}{2+k}\alpha^{(i)}\right)^{1/k} < \beta^{(i)} < \left(\alpha^{(i)}\right)^{1/k}$ for $1 \le i \le r$. The existence of such a $\beta$ follows from the above remark by taking $\eta_i = \frac{1}{2}\left[\left(\frac{2}{2+k}\alpha^{(i)}\right)^{1/k} + \left(\alpha^{(i)}\right)^{1/k}\right]$ for $1 \le i \le r$ and $\varepsilon$ small.

**3. The upper bounds for *g(K, k)*** The upper bounds for $\gamma(s, k)$ implied by the proofs of theorems 2 and 3 are only of interest as an existence proof. However, we can produce a better upper bound by a trivial remark. From corollary 3.1 we have that $g(\boldsymbol{Q}, k)$ is finite, so we can replace the integers $(2k+1)A_j\alpha_j^2$ by sums of at most $g(\boldsymbol{Q},k)$ $k^{\text{th}}$ powers of rational numbers in the final sum. This slight change gives as an upper bound for $\gamma(s, k) \le g(\boldsymbol{Q},k).s.M(s,k)$. The precise value of $g(\boldsymbol{Q}, k)$ is only known for three values of $k$: $g(\boldsymbol{Q}, 2) = 4$, $g(\boldsymbol{Q}, 3) = 3$ and $g(\boldsymbol{Q},$

*4) = 15* (see [8]). For odd *k,* only $w(\boldsymbol{Q}, 3) = 3$ is known. It is obvious that, with the usual notation of the classic Waring problem for $\boldsymbol{Z}$, $g(\boldsymbol{Q}, k) \leq G(k)$. The best known estimate for $G(k)$ (see [18]) is: $G(k) < k(\log k + \log\log k + 2 + O(\log\log k / \log k))$. No general lower bound for $w(\boldsymbol{K}, k)$ is known.

**4. Sums of $k^{th}$ powers in $K(X)$.** If $\boldsymbol{K}$ is a formally real field, we need to know that $w(\boldsymbol{K}(X), 2)$ is finite and that hypothesis (ii) of theorem 3 is satisfied for $\boldsymbol{K}(X)$. This involves showing that certain field quantities are totally positive. In order to simplify the enunciations and proofs of our theorems we need a result of the type: "*$f(X) \in \boldsymbol{K}(X)$ is totally positive if and only if $f(x) \geq 0$ for all $x \in \boldsymbol{K}$.*" Unfortunately such a statement is false in general without some restrictions upon $\boldsymbol{K}$. One such restriction, due to Artin, is that $\boldsymbol{K}$ is a formally real field with precisely one ordering and this ordering is Archimedean. For example, we can take $\boldsymbol{K} = \boldsymbol{Q}$ or $\boldsymbol{R}$ or a finite algebraic extension of $\boldsymbol{Q}$ with precisely one real conjugate field. From now on $\boldsymbol{K}$ will be such a field and we will consider it as a subfield of $\boldsymbol{R}$. With this restriction upon $\boldsymbol{K}$, the hypothesis (ii) of theorem 3 becomes: "If $f(X) \in \boldsymbol{K}(X)$ is such that $f(x) \geq \mu > 0$ for all $x \in \boldsymbol{K}$, there is a rational function $b(X) \in \boldsymbol{K}(X)$ such that $1 < \frac{f(x)}{b(x)^k} < 1 + \frac{2k}{s}$ for all $x \in \boldsymbol{K}$". The existence of such a rational function will be inferred from the following lemma. The proof is, in principle, constructive, since one can use the Bernstein polynomials to construct the function, but it is of no real use for finding a suitable $b(X)$ in practice.

**Lemma 2.** *Let $F(X)$ be a strictly positive definite, everywhere defined continuous real valued function on $\boldsymbol{R}$. Suppose that there exist real numbers $a$, $b$, $\delta$, $C$ and a positive definite rational function $h(X) \in \boldsymbol{Q}(X)$, defined for all $x \in \boldsymbol{R}$, such that:*

$$(1) \qquad 0 < a < \frac{F(x)}{h(x)^k} < b < \infty$$

$$(2) \qquad 0 < \delta \leq F(x)$$

*for all $x$ with $x^2 \geq C > 1$. Then, given $\varepsilon > 0$, there exists $\gamma(X) \in \boldsymbol{Q}(X)$ such that*

$$0 < (a - \varepsilon) < \frac{F(x)}{\gamma(x)^k} < (b + \varepsilon) \quad \text{for all } x \in \boldsymbol{R}.$$

**Proof**. The idea of the proof is simple, we use the Weierstrass approximation theorem to construct a rational function $\gamma(X) \in \boldsymbol{Q}(X)$ which is very close to $h(x)$ for all $x \in \boldsymbol{R}$ with $x^2 > C$ and is sandwiched between $\left(\frac{F(x)}{b+\varepsilon}\right)^{1/k}$ and $\left(\frac{F(x)}{a-\varepsilon}\right)^{1/k}$ for all $x \in \boldsymbol{R}$ with $x^2 \leq C$. However, the details can get a little confusing. We note that hypothesis (2) on $F(X)$ together with the fact that $F(X)$ is strictly positive definite, implies that there is a $\Delta > 0$ such that $F(x) \geq \Delta > 0$ for all $x \in \boldsymbol{R}$.

The following function, where $C > 1$ and $m$ a positive integer, will be used frequently:

$$\alpha_m(X) = \left[1 + \left(\frac{X^2}{C+1}\right)^m\right]^{-1}.$$

The principal properties are:

(a) For each positive integer $m$ and all $x \in \boldsymbol{R}$ we have $0 < \alpha_m(x) \leq 1$.

(b) For any $x$ with $x^2 > C+1$, then $\alpha_m(x) \to 0$ as $m \to \infty$

(c) If $x$ is such that $x^2 = C+1$, then $\alpha_m(x) = ½$ for all $m$.

(d) If $x$ is such that $x^2 < C+1$, then $\alpha_m(x) \to 1$ as $m \to \infty$.

**Step (i).** Given any $\varepsilon > 0$, then if $m > m_0(\varepsilon)$ the function $H_m(X) = \{1 - \alpha_m(X)\} h(X)$

satisfies the inequalities: $\dfrac{F(x)}{b+\varepsilon} < H_m^k(x) < \dfrac{F(x)}{a}$ for all $x \in \boldsymbol{R}$ with $x^2 \geq C+3/2$.

**Proof of Step (i)**. We have $H_m^k(x) = \{1-\alpha_m(x)\}^k h^k(x) \leq h^k(x) < \frac{1}{a}F(x)$ and if we take $\varepsilon_1 > 0$ such that $b(1 - \varepsilon_1)^{-k} < b + \varepsilon$, then by taking $m$ sufficiently large, depending only upon $\varepsilon_1$, we have $0 < \alpha_m(x) < \varepsilon_1$ for all $x \in \boldsymbol{R}$ which satisfy $x^2 > C+3/2$.

It then follows that $H_m^k(x) = \{1-\alpha_m(x)\}^k h^k(x) \geq (1-\varepsilon_1)^k h^k(x) \geq \dfrac{(1-\varepsilon_1)^k}{b} F(x) > \dfrac{F(x)}{b}$ for all $x \in \boldsymbol{R}$ with $x^2 \geq C+3/2$.

We now define a continuous function $G(X)$ on the compact set $x^2 \leq C + 3/2$ as follows:

(1) $G(x) = \dfrac{1}{2}\left(\dfrac{1}{a^{1/k}} + \dfrac{1}{b^{1/k}}\right) F(x)^{1/k}$ for $x \in \boldsymbol{R}$ and $x^2 \leq C$.

(2) $G(x) = h(x)$ *for* $x \in \boldsymbol{R}$ and $C + ½ \leq x^2 \leq C + 3/2$

(3) For $x \in \boldsymbol{R}$ that satisfy $C \leq x^2 \leq C+½$, $G(X)$ can be any continuous function satisfying:

(a) $G(x) = \dfrac{1}{2}\left(\dfrac{1}{a^{1/k}} + \dfrac{1}{b^{1/k}}\right) F(x)^{1/k}$ whenever $x^2 = C$.

(b) $G(x) = h(x)$ whenever $x^2 = C + ½$.

(c) $\left\{\dfrac{F(x)}{b}\right\}^{1/k} < G(x) < \left\{\dfrac{F(x)}{a}\right\}^{1/k}$ for $C < x^2 < C + ½$.

By the Weierstrass polynomial approximation theorem, given any $\varepsilon_2 > 0$ there is a polynomial $P(X)$ with rational coefficients such that $|G(x) - P(x)| < \varepsilon_2$ for all $x \in \boldsymbol{R}$ with $x^2 \leq C + 3/2$. Consider the following rational function, with rational coefficients, defined for all $x \in \boldsymbol{R}$:

$$\gamma_m(X) = \alpha_m(X).P(X) + H_m(X)$$

We will now show that if $m$ is sufficiently large, then $\gamma_m(x)$ satisfies the inequalities of the lemma.

**Step (ii)**. Given $\varepsilon_3 > 0$ there exists $m_3(\varepsilon_3)$ such that for all $m > m_3$ we have $0 < | \alpha_m(x).P(x) | < \varepsilon_3$ for all $x \in \boldsymbol{R}$ with $x^2 \geq C + 3/2$.

**Proof of Step (ii)**. We have $\alpha_m(X)P(X) = P(X)\left(1 + \left(\frac{X^2}{C+1}\right)^m\right)^{-1}$. If the degree of $P(X)$ is $r$ and if $m > r$, then it is a bounded function of $x$ for $x^2 \geq C + 3/2$ and, as $m \to \infty$, this maximum value tends to zero.

**Step (iii)** If $\varepsilon_2$ satisfies the inequalities:

$$0 < \varepsilon_2 < \Delta^{1/k}.min\left[\left\{\left(\frac{1}{a-\varepsilon}\right)^{1/k} - \left(\frac{1}{a}\right)^{1/k}\right\}, \left\{\left(\frac{1}{b+\varepsilon}\right)^{1/k} - \left(\frac{1}{b+2\varepsilon}\right)^{1/k}\right\}\right]$$

Then, if $m > m_2(\varepsilon_2)$, we have for all $x \in \boldsymbol{R}$ with $x^2 \geq C + 3/2$ :

$$\left\{\frac{F(x)}{b+2\varepsilon}\right\}^{1/k} < \gamma_m(x) < \left\{\frac{F(x)}{a-\varepsilon}\right\}^{1/k}$$

**Proof of Step (iii).** If $m > m_3(\varepsilon_2)$, we have by step (i):

$$\gamma_m(x) = (1 - \alpha_m(x) \}.h(x) + \alpha_m(x).P(x) \leq \{ 1 - \alpha_m(x) \}.h(x) + \varepsilon_2$$

$$\leq \left\{\frac{F(x)}{a}\right\}^{1/k} + \varepsilon_2 \leq \left\{\frac{F(x)}{a-\varepsilon}\right\}^{1/k}$$

And if $m > m_3(\varepsilon_2)$ we also have:

$$\left\{\frac{F(x)}{b+2\varepsilon}\right\}^{1/k} + \varepsilon_2 < \left\{\frac{F(x)}{b+\varepsilon}\right\}^{1/k} < \{1-\alpha_m(x)\}h(x)$$

And so $\left\{\frac{F(x)}{b+2\varepsilon}\right\}^{1/k} < \{1-\alpha_m(x)\}h(x) + \alpha_m(x).P(x) = \gamma_m(x)$

**Step (iv)**. If $m > m_1(\varepsilon)$ then for all $x \in \boldsymbol{R}$ with $x^2 \leq C + ½$ we have:

$$\left\{\frac{F(x)}{b+\varepsilon}\right\}^{1/k} < \gamma_m(x) < \left\{\frac{F(x)}{a-\varepsilon}\right\}^{1/k}$$

**Proof of step (iv)**. Let $\varepsilon_3$ satisfy the inequality:

$$0 < 2\varepsilon_3 < min\left[inf\left\{G(x) - \left(\frac{F(x)}{b+\varepsilon}\right)^{1/k}\right\}, inf\left\{\left(\frac{F(x)}{a-\varepsilon}\right)^{1/k} - G(x)\right\}\right],$$

where the inf's are over the $x \in \boldsymbol{R}$ with $x^2 \leq C + ½$. If $m > m_1(\varepsilon_3)$ we have:

$| 1 - \alpha_m(x).h(x) | < \varepsilon_3/4$ and $| 1 - \alpha_m(x).P(x) | < \varepsilon_3/4$

for all $x \in \mathbf{R}$ with $x^2 \leq C + ½$, since $\alpha_m(x) \to 1$ uniformly as $m \to \infty$ and $h(x)$ and $P(x)$ are bounded. Thus we have: $|\gamma_m(x) - P(x)| < \varepsilon_3/2$ and $|\gamma_m(x) - F(x)| < \varepsilon_3$ for all $x \in \mathbf{R}$ with $x^2 \leq C + ½$, Since we have the inequalities:

$$\left\{\frac{F(x)}{b+\varepsilon}\right\}^{1/k} + 2\varepsilon_3 < G(X) < \left\{\frac{F(x)}{a-\varepsilon}\right\}^{1/k} - 2\varepsilon_3$$

we can conclude that for all $x \in \mathbf{R}$ with $x^2 \leq C+½$ we have

$$\left\{\frac{F(x)}{b+\varepsilon}\right\}^{1/k} < \gamma_m(x) < \left\{\frac{F(x)}{a-\varepsilon}\right\}^{1/k}$$

**Step (v)**. If $m > m_3(\varepsilon_2)$, then for all $x \in \mathbf{R}$ which satisfy $C \leq x^2 \leq C + ½$ we have:

$$\left\{\frac{F(x)}{b+\varepsilon}\right\}^{1/k} < \gamma_m(x) < \left\{\frac{F(x)}{a-\varepsilon}\right\}^{1/k}$$

**Proof of Step (v)**. For values of $x$ in the above range we have $G(x) = h(x)$, and so: $|P(x) - h(x)| < \varepsilon_2/2$. Thus, $|\gamma_m(x) - h(x)| = |\gamma_m(x).(P(x) - h(x))| < \varepsilon_2/2$ and by hypothesis

$$\left\{\frac{F(x)}{b}\right\}^{1/k} < h(x) < \left\{\frac{F(x)}{a}\right\}^{1/k}, \text{ so that } \left\{\frac{F(x)}{b+\varepsilon}\right\}^{1/k} < \gamma_m(x) < \left\{\frac{F(x)}{a-\varepsilon}\right\}^{1/k}.$$

To complete the proof of lemma 2, take $m > max\{m_0, m_1, m_2, m_3\}$ and $\gamma(X) = \gamma_m(X)$.

**Theorem 4.** *Let $k > 2$ be a positive integer. If $f(X) \in \mathbf{K}(X)$ is positive definite, then a necessary and sufficient condition that $f(X)$ can be written as a sum of $g(\mathbf{K}(X), k)$ $k^{th}$ powers of totally positive elements of $\mathbf{K}(X)$ is that $2k \mid \partial(f)$.*

**Cor. 4.1.** *If $f(X) \in \mathbf{K}[X]$ is positive definite, $k$ odd and $k \mid \partial(f)$, then $f(X)$ is the sum of at most $g(\mathbf{K}, k)$ $k^{th}$ powers of positive definite rational functions.*

**Proof Cor. 1.** Since $f(X)$ is strictly positive definite $2 \mid \partial(f)$ and $k$ odd, $k \mid \partial(f)$ implies that $2k \mid \partial(f)$.

**Cor. 4.2.** *For any integer $k > 0$, a positive definite polynomial $f(X) \in \mathbf{K}[X]$, with degree divisible by $k$, can be written as a sum of at most $w(\mathbf{K}, k)$ $k^{th}$ powers of $\mathbf{K}(X)$.*

**Proof Cor. 2.** We need only treat the case $k$ even. Suppose $k = 2^r k_1$, where $r > 0$ and $k_1$ is odd. By hypothesis $2.2^{r-1}k_1 \mid \partial(f)$, so from the theorem we can write $f(X)$ as a sum of at most $g(k/2)$ $(k/2)^{th}$ powers of positive definite rational functions, say:

$$f(X) = \sum_{i=1}^{n} a_i^{k/2}(X) \quad \text{with } n \leq g(k/2)$$

Each $a_i(X)$ is positive definite and so can be written as a sum of $s$ squares of rational functions

$$a_i^{k/2} = (b_{1i}^2 + \dots b_{si}^2)^{k/2}$$

We now use the corresponding Hilbert identity of Lemma 1, with *(k+1)* replaced by *k/2* and we get a representation of *f(X)* as a sum of $k^{th}$ powers of rational functions.

**Cor. 4.3.** *The positive definite rational functions f(X) = a(X)/b(X), a(X), b(X)* $\in$ ***K**[X] that can be expressed as sums of at most g(**K(X)**, k)* $k^{th}$ *powers, <u>for every positive integer k</u>, are precisely the set of positive definite rational functions with* $\partial(a) = \partial(b)$.

**Proof Cor. 3.** If $\partial(a) = m$ and $\partial(b) = n$ then $f(X) = a(X)/b(X)$ is a sum of $k^{\text{th}}$ powers iff :

$a(X)b^{k-1}(X)/b^k(X)$ is a sum of $k^{\text{th}}$ powers,

$a(X)b^{k-1}(X)$ is a sum of $k^{\text{th}}$ powers iff $k \mid m + (k-1)n$.

This latter condition holds for all $k$ iff $m = n$.

**Proof of theorem 4.** The condition on the degree of *f(X)* is obviously necessary. There is no loss in generality if we suppose that $f(X) \in \mathbf{K}[X]$ and that *f(X)* is strictly positive definite, i.e. there is a $\mu > 0$ such that $f(x) \geq \mu$ for all $x \in \mathbf{K}$. For $f = p/q$ is a sum of $n$ $k^{\text{th}}$ powers in $\mathbf{K}(X)$ if and only if $p.q^{k-1}$ is a sum of $n$ $k^{\text{th}}$ powers and if $k$ is even or, if $k$ is odd and all the summands are positive definite, then if $f(\vartheta) = 0$, with $\vartheta \in \mathbf{R}$, the irreducible polynomial satisfied by $\vartheta$ is a $k^{\text{th}}$ power factor of *f(X)*.

For the fields $\mathbf{K}$ under consideration we note that $w(\mathbf{K}(X), 2) \leq 5$. To prove the theorem it suffices to show that the modified hypothesis (ii) of theorem 3 is satisfied. That is, we can find a rational function $b(X) \in \mathbf{K}(X)$ such that $1 < \frac{f(x)}{b(x)^k} < 1 + \frac{2k}{s}$ for all $x \in \mathbf{K}$. We use lemma 2 with $\varepsilon$ in the range $0 < \varepsilon < k/(2s + 2k)$. If $f(X) = \alpha X^{2mk} + \dots$, with $\alpha > 0$, we take $\beta \in \mathbf{Q}$ such that $1 - \frac{\varepsilon}{2} < \frac{\alpha s}{(s+k)\beta^k} < 1 + \frac{\varepsilon}{2}$ and $h(X) = (\beta X^{2m} + 1)$. Then, for all $x \in \mathbf{K}$ with $x^2 \geq C(\varepsilon)$ we have: $1 - \varepsilon < \frac{sf(x)}{(s+k)h(x)^k} < 1 + \varepsilon$. From lemma 2 there is a $\gamma(X) \in \mathbf{Q}(X)$ such that $1 - 2\varepsilon < \frac{sf(x)}{(s+k)\gamma(x)^k} < 1 + 2\varepsilon$ and since $\varepsilon < k/(2s + 2k)$ we have:

$$1 < (1 - 2\varepsilon)(s + k)/s \text{ and } (1 + 2\varepsilon)(s + k)/s < 1 + k/s.$$

**<u>5. Upper bounds for *w(K(X), k)* and *g(K(X), k)*</u>.** There is little precise information available.

For $k = 2$, $w(\boldsymbol{K}(X), 2) = s(\boldsymbol{K}) + 1$ (see [10] when $\boldsymbol{K}$ is an algebraic number field and it is classic when $\boldsymbol{K} = \boldsymbol{R}$). Choi *et al.* [4] showed that $g(\boldsymbol{R}(X), 2) = 2$. We extend this result in theorem 5 by showing that $w(\boldsymbol{K}(X), 2) = g(\boldsymbol{K}(X), 2)$.

For $k = 3$, if $\boldsymbol{F}$ is any field not of characteristic $3$, $w(\boldsymbol{F}, 3) \leq 3$ and $g(\boldsymbol{K}(X), 3) \leq 3$. These follow from a classic identity due to Richmond (see [8], Notes to Chapter13 and Manin [11]): If $r, s \in \boldsymbol{F}$ and $s \neq 0$, then

$$r = \left\{\frac{s(1+t^3)}{3(1-t+t^2)}\right\}^3 + \left\{\frac{s(3t-1-t^3)}{3(1-t+t^2)}\right\}^3 + \left\{\frac{s(3t-3t^2)}{3(1-t+t^2)}\right\}^3, \quad \text{where } t = 3r/s^3.$$

This gives $w(\boldsymbol{F}, 3) \leq 3$. To deduce that $g(\boldsymbol{K}(X), 3) \leq 3$, we note that $(1 - t + t^2) = (t - ½)^2 + ¾$ is totally positive, so we must show that if $6 \mid \partial r$ and $r(x) > 0$ for all $x \in \boldsymbol{K}$, then $s(X) \in \boldsymbol{K}(X)$ can be chosen so that $s(x) > 0$, $(t(x) - t^2(x)) > 0$ and $3t(x) - 1 - t^3(x) > 0$ for all $x \in \boldsymbol{K}$. A sufficient condition is that $½ < t(x) < 1$ for all $x \in \boldsymbol{K}$, namely $3r < s^3 < 6r$. If $r(X) \in \boldsymbol{K}[X]$ is strictly positive definite, of degree $6m$ and with leading coefficient $a$, then if $h(X) = bX^{2m} + 1$, where $b$ is such that $3a < b^3 < 6a$, we have $1/6 < r(x)/h(x)^3 < 1/3$ for all $x \geq C(f)$. We now use lemma 2 to construct $s(X)$ to conclude that $g(\boldsymbol{K}(X), 3) \leq 3$.
For $k = 4$, Choi *et al.* [4] use the fact that $g(\boldsymbol{R}(X), 2) = 2$ to show that $w(\boldsymbol{R}(X), 4) \leq 6$.

As a final application of lemma 2 we prove the following theorem.
**Theorem 5**. *(i) If $4 \mid \partial f$ and for some $n$, $f(X) = a_1^2 + \ldots + a_n^2$ where $a_i(X) \in \boldsymbol{K}(X)$, then $f(X)$ can be represented as $f(X) = g_1^2 + \ldots + g_n^2$ where, $1 \leq i \leq n$, $g_i(X) \in \boldsymbol{K}(X)$ and $g_i(x) \geq 0$ for all $x \in \boldsymbol{K}$.*
*(ii) If $4 \mid \partial f$, $f(x) \geq \mu_0 > 0$ for all $x \in \boldsymbol{K}$ and for some $n$ $f(X) = a_1^2 + \ldots + a_n^2$ where $a_i(X) \in \boldsymbol{K}(X)$, then $f(X)$ can be represented as : $f(X) = g_1^2 + \ldots + g_n^2$ where, $1 \leq i \leq n$, $g_i(X) \in \boldsymbol{K}(X)$ and $g_i(x) \geq \mu_i > 0$ for all $x \in \boldsymbol{K}$.*

**Cor. 1** *Either $w(\boldsymbol{K}(X), 2) = \infty$ or $w(\boldsymbol{K}(X), 2) = g(\boldsymbol{K}(X), 2)$.*
**Cor. 2.** *If $w(\boldsymbol{K}(X), 2) = 2$, then $w(\boldsymbol{K}(X), 4) \leq 20$. If $\sqrt{2} \in \boldsymbol{K}$, then $w(\boldsymbol{K}(X), 4) \leq 8$. If $\sqrt{3} \in \boldsymbol{K}$, then $w(\boldsymbol{K}(X), 4) \leq 6$.*
**Proof of Cor. 2**. If $f(X) \in \boldsymbol{K}[X]$ is strictly positive definite and $4 \mid \partial f$, then $f(X)/6$ can be represented as $f(X)/6 = g_1^2 + g_2^2$ with $g_1$ and $g_2$ strictly positive definite. We then have $g_1 = a^2 + b^2$ and $g_2 = c^2 + d^2$. The Hilbert identity $6(a^2 + b^2)^2 = (a + b)^4 + (a - b)^4 + 4a^4 + 4b^4$ applied twice gives the first result. If $\sqrt{2} \in \boldsymbol{K}$, then $4 = (\sqrt{2})^4$ and we have the second result.
If $\sqrt{3} \in \boldsymbol{K}$, then we use the identity $18(a^2 + b^2)^2 = (a + \sqrt{3}b)^4 + (a - \sqrt{3}b)^4 + (2b)^4$.

**Proof of theorem 5**. We first show that the truth of (ii) implies the truth of (i).

Suppose that $f(x) \geq 0$ for all $x \in \boldsymbol{K}$ and that $f(X) \in P(\boldsymbol{K}(X), 4)$. We have

$f(X) = a_1^4(X) + .... + a_r^4(X)$ and if $\gamma \in \boldsymbol{R}$ is such that $f(\gamma) = 0$, then $a_i(\gamma) = 0$ for $i = 1, ..., n$ and this is true for all the conjugates of $\gamma$. Hence the irreducible polynomial $p(X)$ over $\boldsymbol{Q}$ defining $\gamma$ must divide each $a_i(X)$, which implies the $p^4(X)$ divides $f(X)$. Repeating the process we have $f(X) = P^4(X)F_1(X)$, where $F_1(X) \geq \mu > 0$ for all $x \in \boldsymbol{K}$.

If $f(X) = b_1^2(X) + .... + b_n^2(X)$, then $P^2(X)$ divides each $b_i(X)$ and we have a representation of $F_1(X)$ as a sum of $n$ squares. Since $4 \mid \partial F_1$, by (ii), we have $F_1(X) = g_1^2(X) + .... + g_n^2(X)$ where $g_i(x) \geq \mu_i > 0$ for all $x \in \boldsymbol{K}$. We then have the representation:

$$f(X) = P^4(X)F_1(X) = (P^2(X)g_1(X))^2 + .... + (P^2(X)g_n(X))^2,$$

where, for $1 \leq i \leq n$, $P^2(x)g_i(x) \geq 0$ for all $x \in \boldsymbol{K}$.

We shall prove (ii) by induction on $n$. If $n = 1$, then $f(X) = a_1^2(X)$ the result is immediate.

The induction hypothesis is:

"(ii) holds for all $f(X) \in P(\boldsymbol{K}(X), 4)$ and such that $f(X)$ can be written as a sum of at most $(n-1)$ squares in $\boldsymbol{K}(X)$".

The above remarks imply the following induction hypothesis:

"(i) holds for all $f(X) \in P(\boldsymbol{K}(X), 4)$ and such that $f(X)$ can be written as a sum of at most $(n-1)$ squares in $\boldsymbol{K}(X)$". This will be used in the proof of lemma 5.

In order to prove the 'case n' of the induction hypothesis we need the following four lemmas.

**Lemma 3**. *If **K** is a field, not of characteristic 2, and* $f = a_1^2 + a_2^2 + ... + a_n^2$, $f$, $a_i \in \boldsymbol{K}$ *for* $1 \leq i \leq n$, *then the general solution of* $fU_0^2 = U_1^2 + ... + U_n^2$ *is given by:*

$$u_0 = \sum_{j=1}^{n} T_j^2 \quad u_i = 2T_i(\sum_{j=1}^{n} a_j T_j) - a_i \sum_{j=1}^{n} T_j^2 \quad for\ 1 \leq i \leq n,$$

*where* $T_j \in \boldsymbol{K}$ *for* $1 \leq i \leq n$.

**Proof of lemma 3.** If $\underline{\boldsymbol{a}} = (a_1, ..., a_n)$, where $f = a_1^2 + a_2^2 + ... + a_n^2$, $\underline{\boldsymbol{t}} = (t_1, ..., t_n) \in \boldsymbol{K}^n$, then the line joining $\underline{\boldsymbol{a}}$ to $\underline{\boldsymbol{t}}$ intersects the quadric of $fU_0^2 = U_1^2 + ... + U_n^2$ in a second $K$-rational point. If we write $T_i = t_i - a_i$ for $i = 1, ..., n$, then the coordinates of this $K$-rational point are given by the above formulae.

**Lemma 4.** *If* $f(X) \in \boldsymbol{K}[X]$ *is strictly positive definite,* $4 \mid \partial f$ *and* $f = a_1^2 + ... + a_n^2$, $a_i \in \boldsymbol{K}(X)$, *then there exists* $g_0(X), g_1(X), ..., g_n(X) \in \boldsymbol{K}[X]$ *such that* $f.g_0^2 = g_1^2 + ... + g_n^2$, $\partial g_1 = ... = \partial g_n$, *and* $g_0(X)$ *is strictly positive definite.*

**Proof of lemma 4.** By hypothesis, $f = a_1^2 + ... + a_n^2$, $a_i \in \boldsymbol{K}(X)$. By Cassels' theorem [2], there are $b_1, ..., b_n \in \boldsymbol{K}[X]$ such that $f = b_1^2 + ... + b_n^2$. Since $\partial f = 4m$, at least one of the $b_i$, say $b_1$, has $\partial b_1 = 2m$ and $\partial b_i \leq 2m$ for $i = 2, ..., n$.

In lemma 3, choose $T_1$ to be a polynomial of degree $d > 1$ and, for $j > 1$, $T_j$ a polynomial of degree $< d$. We then have, for $i = 1, ..., n$, $\partial h_i = 2m + d$.

**Lemma 5.** *If $f(X) \in \mathbf{K}[X]$ is strictly positive definite, $4 \mid \partial f$ and $f.a_0^2 = a_1^2 + ... + a_n^2$, $a_i \in \mathbf{K}[X]$, $a_0(X)$ strictly positive definite and $\partial a_1 = ... = \partial a_n$, then there exists $g_0(X), g_1(X), ..., g_n(X) \in \mathbf{K}[X]$ such that $f.g_0^2 = g_1^2 + ... + g_n^2$, $\partial g_1 = ... = \partial g_n$, $g_0(X)$ is strictly positive definite and $g_1(X)$ is positive definite.*

**Proof of lemma 5**. From the representation $f.a_0^2 = a_1^2 + ... + a_n^2$ we have $F(X) = a_0^2 f(X) - a_n^2 = a_1^2 + ... + a_{n-1}^2$ and since $4 \mid \partial F(X)$, by the induction hypothesis we have $c_0^2 F = c_1^2 + ... + c_{n-1}^2$, where $c_0(x) \geq \mu > 0$ and for $i = 1, ..., n$, $c_i(x) \geq 0$ for all $x \in \mathbf{K}$. Hence:

$(a_0c_0)^2 f = c_1^2 + ... + c_{n-1}^2 + (c_0a_n)^2$ and we can take $g_0 = a_0c_0$ and $g_1 = c_1$.

**Lemma 6.** *If $f(X) \in \mathcal{P}(\mathbf{K}[X], 4)$ is strictly positive definite and $f.a_0^2 = a_1^2 + ... + a_n^2$, $a_i \in \mathbf{K}[X]$, $a_0(X)$ strictly positive definite, $a_1(X)$ positive definite and $\partial a_1 = ... = \partial a_n$, then there exists $g_0(X), g_1(X), ..., g_n(X) \in \mathbf{K}[X]$ such that $f.g_0^2 = g_1^2 + ... + g_n^2$, $\partial g_1 = ... = \partial g_n$ and $g_0(X)$, $g_1(X)$ are both strictly positive definite.*

**Proof of lemma 6.** By hypothesis $f = a_1^2 + ... + a_n^2$, where $a_1(x) \geq 0$ for all $x \in \mathbf{K}$ and $\partial a_1 = \partial a_n = 2m$. The $a_i(X)$ cannot have any common real zeros since $f(x) \geq \mu > 0$ for all $x \in \mathbf{K}$. In the general solution of $fU_0^2 = U_1^2 + ... + U_n^2$ we choose $T_1 = 1$ and for $j = 2, ..., n$, we choose

$T_j = a_j/(AX^{2m} + C)$, where $A$, $C$ are positive rational numbers chosen so that $\frac{a_2^2 + .... + a_n^2}{\left(Ax^{2m} + C\right)^2} \leq \frac{1}{2}$ for all $x \in \mathbf{K}$ and we write $b(x) = 1 - \frac{a_2^2 + .... + a_n^2}{\left(Ax^{2m} + C\right)^2} \geq \frac{1}{2}$

From lemma 3 we have:

$g_0(x) = 1 + T_2^2 + .... + T_n^2 \geq 1$ for all $x \in \mathbf{K}$. $g_1(x) = a_1(x)b(x) + \frac{2}{Ax^{2m} + C}\sum_{j=2}^{n} a_j^2(x)$

Thus $g_1(x)$ is the sum of $n$ positive terms and so $g_1(x) \geq 0$ for all $x \in \mathbf{K}$. If $\gamma \in \mathbf{R}$ is such that $g_1(\gamma) = 0$ we must have $a_1(\gamma) = ... = a_n(\gamma) = 0$. We have seen that this is impossible, hence, $g_1(x) \geq \mu_1$ for all $x \in \mathbf{K}$.

We are now in a position to prove the induction step from '*n-1*' to '*n*'.

Suppose that $f(X) \in P(\boldsymbol{K}(X), 4), f(x) \geq \mu > 0$ for all $x \in \boldsymbol{K}$ and $f(X) = a_1^2(X) + ... + a_n^2(X)$. By renumbering the functions $a_i(X)$ and using lemmas 4, 5 and 6 if necessary, we deduce that we have a representation:

$b_0^2(X)f(X) = b_1^2(X) + ... + b_n^2(X)$, *where* $b_0(x) \geq \mu_0 > 0$ *and* $b_1(x) \geq \mu_1 > 0$ *for all* $x \in \boldsymbol{K}$ and $\partial b_1 = ... = \partial b_n$. We shall deduce that there is a representation:

$g_0^2(X)f(X) = g_1^2(X) + ... + g_n^2(X)$, *where, for* $i = 0, 1, ..., n$, $g_i(x) \geq \mu_i > 0$ *for all* $x \in \boldsymbol{K}$.

Let $\lambda(X) = b_1^2(X)/2f(X)$, then, since $2\partial b_1 = \partial f$, $\lambda(x) \geq \lambda_0 > 0$ for all $x \in \boldsymbol{K}$. From lemma 2, for any $\lambda_0 > 0$, there exists $\phi(X) \in \boldsymbol{Q}(X)$ such that $f(x) < \phi^4(x) < (1+\lambda_0)f(x)$ for all $x \in \boldsymbol{K}$.

We use the formulae of lemma 3 to construct a new representation $h_0^2(X)f(X) = h_1^2(X) + ... + h_n^2(X)$, by taking $T_1 = 1$, $T_2 = (b_2(X) + \phi^2(X))/b_1(X)$ *and* $T_j = b_j(X)/b_1(X)$ *for* $j = 3, ..., n$.

We have:

$h_0(X) = 1 + T_1^2 + .... + T_n^2 \geq 1$ for all $x \in \boldsymbol{K}$.

$$h_1 = 2b_2T_2 + 2\sum_{j=3}^{n} b_jT_j - b_1(1 - T_2^2 - .... - T_n^2) = \frac{\phi^4 - f}{b_1} \geq \mu_1 > 0,$$

$$\frac{h_1^2}{h_0^2} \leq \frac{(\phi^4 - f)^2}{h_0^2 b_1^2} \leq \frac{\lambda_0^2 f^2}{h_0^2 b_1^2} \leq \frac{b_1^4 f^2}{4f^2 h_0^2 b_1^2} = \frac{b_1^2}{4h_0^2} \leq \frac{b_1^2}{4} \leq \frac{f}{4}.$$

Thus $f_1 = f - \dfrac{h_1^2}{h_0^2} = \left(\dfrac{h_2}{h_0}\right)^2 + .... + \left(\dfrac{h_n}{h_0}\right)^2 \geq \dfrac{3}{4} f \geq \dfrac{3\mu}{4} > 0$. Since $4 \mid \partial f_1$ and $f_1$ is strictly positive definite, by theorem 4, $f_1 \in P(\boldsymbol{K}(X), 4)$, also $f_1$ can be represented as a sum of $(n-1)$ squares, so by our induction hypothesis there exist $g_0, g_1, ..., g_{n-1}$ such that $g_0^2(X)f_1(X) = g_1^2(X) + ... + g_{n-1}^2(X)$, where, for $i = 0, 1, ..., n-1$, $g_i(x) \geq \mu_i > 0$ *for all* $x \in \boldsymbol{K}$. Hence: $(g_0h_0)^2 f = (g_0h_1)^2 + g_1^2 + .... + g_{n-1}^2$ and we have our desired representation for $f$. This completes the induction step and the proof of theorem 5.

**References**

[1] E. Becker, *Summen n-ter Potenzen in Körpern*, J. reine angew. Math. 307/308 (1979), 8-30.

[2] J.W.S. Cassels, *On the representation of rational functions as sums of squares*, Acta Arith. 9 (1964), 79-82.

[3] J.W.S. Cassels, W.J. Ellison, A. Pfister, *On sums of squares and on elliptic curves over function fields*, J. Number Theory, 3 (1971), 125-149.

[4] M.D. Choi, T.Y. Lam, A. Prestel, B. Reznick, *Sums of* $2m^{th}$ *powers of rational functions in one variable over real closed fields,* Math. Zeit. 221, (1996), 93-112.

[5] J.L Colliot-Thénène, U. Jannsen, *Sommes de carrés dans les corps de fonctions*, C.R.A.S. 312 (1991), 759-762. (Now that the Milnor conjectures have been proved.)

[6] W.J. Ellison, *Waring's Problem*, Amer. Math. Monthly, 78, N° 1 (1971), 10-36.

[7] W.J. Ellison, *Waring's and Hilbert's 17$^{th}$ problem*, Ph.D. Dissertation, 1970, University of Cambridge, N° D1537/71.

[8] G.H. Hardy and E.M. Wright, *An Introduction to the Theory of Numbers*, Oxford, 1960.

[9] D.W. Hoffmann, *Pythagoras numbers of fields*, J. Amer. Math. Soc. 12, N°2 (1999), 839-848.

[10] J.S. Hsia, R.P. Johnson, *On the representation in sums of squares for definite functions in one variable over an algebraic number field*, Amer. J. of Math. 96 N°3 (1974), 448-453.

[11] Yu. Manin, *Cubic Forms: Algebra, Geometry, Arithmetic*, Elsevier Science Publishing, 1986.

[12] Yu. V. Nesterenko, *On Waring's problem (elementary methods)*, Journal of Mathematical Sciences, Vol. 137, No. 2, 2006. Translated from Zapiski Nauchnykh Seminarov POMI, Vol. 322, 2005, 149–175.

[13] A. Oppenheim, *Hilbert's proof of Waring's theorem*, Messenger of Math. 58 (1929), 153-158.

[14] A. Pfister, *Zur Darstellung definiter Funktionen als Summe von Quadraten*, Invent. Math. 4 (1967), 229-237.

[15] A. Pfister, *Quadratic Forms with Applications to Algebraic Geometry and Topology*, London Math. Soc. Lecture Note series N° 217, Cambridge University Press, 1995.

[16] Y. Pourchet, *Sur la représentation en somme de carrés des polynômes à une indéterminée sur un corps de nombres algébriques*, Acta Arith. 19 (1971), 89-104.

[17] E. Stridsberg, *Sur la démonstration de M. Hilbert du théorème de Waring*, Math. Ann., 74 (1913), 271-274.

[18] R.C. Vaughn and T.D. Wooley, *Waring's problem: A survey*, in Number theory for the millennium, 2002.

W.J. Ellison
15 allée de la Borde
33450 St. Sulpice et Cameyrac, FRANCE
e-mail : w.ellison@wanadoo.fr